\documentclass[reqno,11pt]{amsart}

\usepackage{amsmath, amssymb, amsfonts, amsthm}  % AMS 数学包
\usepackage{mathrsfs}    % 花体字母
\usepackage{bbm}         % 黑板粗体 \mathbbm
\usepackage{dsfont}      % \mathds
\usepackage{latexsym}    % 经典符号
\usepackage{euscript}    % \EuScript
\usepackage{hyperref}
\usepackage{geometry}
\allowdisplaybreaks      % 允许公式跨页断开

\usepackage{xcolor}      % 颜色支持
\usepackage{enumitem}    % 可定制列表
\usepackage{enumerate}   % enumerate 扩展
\makeatletter

\newcommand{\Rmnum}[1]{\expandafter\@slowromancap\romannumeral #1@}
\makeatother

\def\bos{\boldsymbol}

\numberwithin{equation}{section}  % 公式按 section 编号

\newtheorem{theorem}{Theorem}[section]
\newtheorem{lemma}[theorem]{Lemma}
\newtheorem{corollary}[theorem]{Corollary}
\newtheorem{definition}[theorem]{Definition}
\newtheorem{proposition}[theorem]{Proposition}
\newtheorem{remark}[theorem]{Remark}

\begin{document}
\title[Hausdorff dimension of the Cartesian product of exact approximation set]{Hausdorff dimension of the Cartesian product of exact approximation set in $\beta$-expansions}

\author[W. Cheng]{Wanjin Cheng}
\address[Wanjin Cheng]{School of Mathematics, South China University of Technology, Guangzhou, 510640, China}
\email{chengwj0227@163.com}
\author[X. Zhang]{Xinyun Zhang$^{\star}$}
\address[Xinyun Zhang]{School of Mathematics and Information Science, Nanchang Hangkong University, Nanchang, 330063, China}
\email[corresponding author]{xinyunzhang@nchu.edu.cn}


\begin{abstract}
In this paper, we study the metrical theory of Cartesian products of exact approximation sets in $\beta$-expansions.
More precisely, for an integer $d \ge 2$ and real numbers $\beta_i > 1$ $(1 \le i \le d)$,
we consider the set of points $x_i \in [0,1)$ is approximable
by its convergents in the $\beta_i$-expansion to order $\psi_i$, but not to any better order.
For any non-increasing functions $\psi_i$, we determine the Hausdorff dimension of the Cartesian product of these sets.
\end{abstract}
\maketitle

\section{introduction}
Metric Diophantine approximation, a central topic in number theory, concerns how well real numbers can be approximated by rationals, with foundational contributions due to Khintchine and Jarn\'ik.

%%From a global, measure-dimension perspective, this line of inquiry belongs to the theory of metric Diophantine approximation. This is originated from the work of Khintchine and Jarn\'ik.
%\textcolor{blue}{
%Khintchine~\cite{Khintchine1936} investigated the measure theory of the
%$\psi$-approximable set
%\[
%W(\psi) := \Bigl\{ x \in \mathbb{R} : \bigl|x - \tfrac{p}{q}\bigr| < \psi(q)
%\ \text{for infinitely many } (p,q)\in\mathbb{Z}\times\mathbb{N} \Bigr\},
%\]
%where $\psi : \mathbb{N} \to \mathbb{R}^+$ is a decreasing function such that $x \mapsto x^2 \psi(x)$ is non-increasing. He established a full/zero dichotomy law:
%\begin{equation*}
%\[\mathcal{L}\big(\mathcal{W}(\psi)\cap[0,1]\big)=
%\begin{cases}
%& 0, \quad\quad\quad\quad\;\text{if}~
%\sum\limits_{q=1}^{\infty}q\cdot\psi(q)
%<+\infty,\\
%& 1,  \quad\quad\quad\quad\;\text{if}~
%\sum\limits_{q=1}^{\infty}q\cdot\psi(q)
%=+\infty,
%\end{cases}\]
%\end{equation*}
%where $\mathcal{L}$ denotes the Lebesgue measure. }
%This theory originates from the work of Khintchine and Jarn\'ik.
For a decreasing function $\psi : \mathbb{N} \to \mathbb{R}^+$ with $x \mapsto x^2 \psi(x)$ non-increasing,
Khintchine~\cite{Khinchin24} established a zero--one law for the Lebesgue measure of the set of $\psi$-well approximable numbers:
\[
W(\psi) := \Bigl\{ x \in \mathbb{R} : \bigl|x - \tfrac{p}{q}\bigr| < \psi(q)
\ \text{for infinitely many } (p,q)\in\mathbb{Z}\times\mathbb{N} \Bigr\},
\]
that is, $\mathcal{L}(W(\psi)) = 0 \ \text{or}\ 1$ according as
$\sum_{q \ge 1} q \psi(q)$ converges or diverges.
Jarn\'ik~\cite{Jarnik29} and Besicovitch~\cite{Besicovitch34} independently showed that for $\tau \ge 2$,
$\dim_{\mathcal{H}} W\!\left(x \mapsto x^{-\tau}\right) = \frac{2}{\tau}$.
Here and throughout, $\mathcal{L}$ and $\dim_{\mathcal{H}}$ denote the Lebesgue measure and the Hausdorff dimension, respectively.
%\textcolor{blue}{ To study the exceptional sets in Diophantine approximation (sets with Lebesgue measure zero),  Jarn\'ik pioneered the application of fractal measure theory. He~\cite{Jarnik1929} and Besicovitch~\cite{Besicovitch1934} independently showed that
%%Jarn\'ik~\cite{Jarnik1929} and independently Besicovitch~\cite{Besicovitch1934}
%$$\dim_\mathcal{H} W\!\left(x \mapsto x^{-\tau}\right) = \frac{2}{\tau}\ \  \text{for\ all}\ \tau \geq 2.$$
%Here and hereafter, $\dim_\mathcal{H}$ denotes the Hausdorff dimension.}
Subsequently,  for a general dimension function $f$, Jarn\'ik~\cite{Jarnik31} determined the Hausdorff $\mathcal{H}^f$-measure of the set $W(\psi)$. %For a modern version of these two classic results, we refer the reader to Beresnevich, Dickinson, and Velani~\cite{BDV1, BDV2}.

Jarn\'ik \cite{Jarnik31} also introduced the set of points with exact order of approximation, defined by
\begin{equation}\label{eq:exact}
\mathrm{Exact}(\psi) := W(\psi) \setminus \bigcup_{0 < c < 1} W(c\psi).
\end{equation}
More precisely, $\mathrm{Exact}(\psi)$ consists of the points that belong to $W(\psi)$ but do not lie in any of the sets $W(c\psi)$ for any $0 < c < 1$.
Jarn\'ik showed that $\mathrm{Exact}(\psi)$ is non-empty whenever $\psi(q) = o(q^{-2})$.

Regarding the metric theory of $\mathrm{Exact}(\psi)$, Bugeaud \cite{Bugeaud03} proved that if $\psi$ is non-increasing and satisfies $\psi(q) = o(q^{-2})$, then
\[
\dim_{\mathcal{H}} \mathrm{Exact}(\psi)
= \dim_{\mathcal{H}} W(\psi)
= \frac{2}{\lambda},
\quad
\lambda = \liminf_{n \to +\infty} \frac{-\log \psi(n)}{\log n}.
\]

%where $\lambda:=\liminf\limits_{n\rightarrow+\infty}\frac{-\log\psi(n)}{\log n}$.}
%While the theory of exact approximation in the classical setting has been well developed, much less is known in the context of $\beta$-expansions.
In this paper, we investigate the Cartesian products of exact
approximation sets  under the framework of $\beta$-expansions.
\subsection{Approximation in $\beta$-expansions.}
Diophantine approximation problems have also been studied widely in the setting of
$\beta$-expansions. We begin with a brief introduction to $\beta$-expansions.

Let $\beta > 1$ be a real number, and $T_\beta : [0,1] \to [0,1]$ be the
$\beta$-transformation defined as
$$
T_\beta(x) = \beta x - \lfloor \beta x \rfloor,
$$
where $\lfloor \cdot \rfloor$ represents the integer part of a real number. Then every $x \in [0,1]$ can be uniquely expanded
into a finite or an infinite series
\begin{equation}\label{eq1}
x=\frac{\varepsilon_{1}(x)}{\beta}+\frac{\varepsilon_{2}(x)}{\beta^2}+\cdots+
\frac{\varepsilon_{n}(x)+T_{\beta}^{n}(x)}{\beta^{n}}
=\sum_{n=1}^{+\infty}\frac{\varepsilon_{n}(x)}{\beta^{n}},
\end{equation}
where $\varepsilon_{1}(x)=\lfloor\beta x\rfloor$ and $\varepsilon_{n+1}(x)=\varepsilon_{1}\big(T_{\beta}^{n}(x)\big)$ ($n\in\mathbb{N}$), are called the digits of the $\beta$-expansion of $x.$ %Sometimes $x$ can be identified with its \emph{digit} sequence
%$$\varepsilon(x):=\big(\varepsilon_{1}(x),\ldots,\varepsilon_{n}(x),\ldots\big).$$
%If there exists some $n_{0}\in\mathbb{N}$ such that $\epsilon_{n}(x)=0$ for all $n\geq n_{0},$ we say that the $\beta$-expansion of $x$ is \emph{finite}. Otherwise, it is said to be \emph{infinite}.}

For a real number $x\in[0, 1],$ we call the partial sums of the series
$$\omega_{n}(x)=\frac{\varepsilon_{1}(x)}{\beta}+\frac{\varepsilon_{2}(x)}{\beta^{2}}+
\cdots+\frac{\varepsilon_{n}(x)}{\beta^{n}}$$
the \emph{convergents} of the $\beta$-expansion of $x$. Since $T_\beta^n(x) \in [0,1)$, it follows from \eqref{eq1} that
\[
\bigl| x - \omega_n(x) \bigr| < \beta^{-n},
\quad \forall n \in \mathbb{N}.
\]

%\textcolor{blue}{In analogy to classical Diophantine approximaiton, Philipp \cite{Philipp} determined the Lebesgue measure of the set
%$$W_{\beta}(\psi):=\Big\{x\in \big[0,1\big]: \big|T_{\beta}^n(x)-y\big|<\psi(n) \ {\text{for infinitely many}}\ n\in \mathbb{N}\Big\},$$ where $y \in[0, 1]$.  Shen and Wang \cite{ShenWang13} determined the Hausdoff dimension of $W_{\beta}(\psi)$.}
%An extensively studied topic in $\beta$-expansions is the shrinking target problems
%in a dynamical system, which began with the pioneering work of Hill and Velani~\cite{HillVelani},
%concerning the Diophantine properties of the orbits.
%As far as $\beta$-expansion is concerned, Philipp~\cite{Philipp} determined the Lebesgue measure
%of the set
%\[
%W(\psi) := \Bigl\{ x \in [0,1] : \, |T_\beta^n(x) - y| < \psi(n)
%\ \text{for infinitely many } n \in \mathbb{N} \Bigr\},
%\]
%where $y \in [0,1]$. The Hausdorff dimension of $W(\psi)$ was given by Shen and Wang~\cite{ShenWang}
%as follows.

%\begin{theorem}[{\cite[Theorem~1.1]{ShenWang}}]
%Let $\psi$ be a positive function defined on $\mathbb{N}$. For any $\beta > 1$,
%\[
%\dim_{H} W(\psi) = \frac{1}{1+b},
%\]
%with
%\[
%b = \liminf_{n \to \infty} \frac{- \log_\beta \psi(n)}{n},
%\]
%where $\dim_{H}$ denotes the Hausdorff dimension.
%\end{theorem}

%% More precisely, for a real number $x \in [0,1]$, we call
%%\[\omega_n(x) = \frac{\epsilon_1(x)}{\beta}+ \frac{\epsilon_2(x)}{\beta^2}+ \cdots+ \frac{\epsilon_n(x)}{\beta^n},\qquad n \in \mathbb{N},\]be the \emph{convergents} of the $\beta$-expansion of $x$.

A natural question is: how much faster can $\omega_n(x)$ approach to $x$ than the trivial bound $\beta^{-n}$?
Fang et. al. \cite{Fang20} investigated this problem and proved that, for any $\beta > 1$,
$$
\lim_{n \to +\infty} \frac{1}{n} \log_\beta \bigl| x - \omega_n(x) \bigr|
= -1 \quad \text{for}\ \mathcal{L}\ \text{almost every}\ x\in[0,1).$$
Moreover, if $\varphi$ is non-decreasing and satisfies
$\eta := \liminf\limits_{n \to +\infty} \frac{\varphi(n)}{n} \geq 1,$
then
\[
\dim_\mathcal{H} \Bigl\{ x \in [0,1] :
\liminf_{n \to +\infty} \frac{1}{\varphi(n)}
\log_\beta \bigl| x - \omega_n(x) \bigr| = -1 \Bigr\}
= \frac{1}{\eta}.
\]

Recently, Zhang and Zhong~\cite{Zhang24} generalized the problem of exact Diophantine approximation in $\beta$-expansions. They defined
\[W_\beta(\psi) := \Bigl\{ x \in [0,1] :\bigl|x - \omega_n(x)\bigr| < \frac{\psi(n)}{\beta^n}\ \text{for infinitely many } n \in \mathbb{N} \Bigr\},\]
and analogous to \eqref{eq:exact},
\[
E_\beta(\psi) := W_{\beta}(\psi) \setminus \bigcup_{0<c<1} W(c\psi).
\]
They determined the Hausdorff dimension of $E_\beta(\psi)$.

\begin{theorem}[\cite{Zhang24}]\label{thzhang}
Let $\beta > 1$ be a real number and let $\psi : \mathbb{N} \to \mathbb{R}_{+}$
be a positive and non-increasing function with $\psi(n) \to 0$ as $n \to +\infty$.
Then
\[
\dim_\mathcal{H} E_\beta(\psi) = \frac{1}{1+\alpha},\quad
\alpha= \liminf_{n \to +\infty} \frac{-\log_\beta \psi(n)}{n}.
\]
\end{theorem}

\subsection{The product of approximation sets}In this section, we turn to the Cartesian product of approximation sets and study their metric properties. The foundational results concerning the fractal dimensions of Cartesian product sets were obtained by Marstrand \cite{Marstrand54} and Tricot \cite{Tricot82}.
\begin{theorem}[\cite{Marstrand54}, \cite{Tricot82}]
Suppose that $A \subset \mathbb{R}^d$ and $B \subset \mathbb{R}^t$ are two Borel measurable sets. Denote by $\dim_\mathcal{P}$ the packing dimension. Then
\[
\dim_\mathcal{H} A + \dim_\mathcal{H} B \leq \dim_\mathcal{H} (A \times B) \leq \dim_\mathcal{H} A + \dim_\mathcal{P} B.
\]
\end{theorem}

The study of the Hausdorff dimension of Cartesian products of approximation sets dates back to Schleischitz \cite{Schleischitz23}, motivated by Erd\H{o}s' work~\cite{Erdos62}.
It is well known that the set of Liouville numbers is rather small in the sense that it has Hausdorff dimension $0$.
Nevertheless, Erd\H{o}s~\cite{Erdos62} proved that every real number can be written as the sum of two Liouville numbers. Consider the Lipschitz map $f: L \times L \to \mathbb{R}$ defined by $f(x,y) = x+y$,
by the property of Hausdorff dimension, it follows that
\[
\dim_{\mathcal{H}} (L \times L) = 1,
\]
where
\[
L := \bigcap_{v > 1} \Bigl\{ x \in \mathbb{R} : \| qx \| < q^{-v} \text{ for infinitely many } q \in \mathbb{N} \Bigr\},
\]
and $\|\cdot\|$ denotes the distance to the nearest integer.

%In 1962, Erd\"{o}s \cite{erd} demonstrated that the Hausdorff dimension of the Cartesian product of the set of Liouville numbers satisfies

Following this insight, Schleischitz \cite{Schleischitz23} obtained bounds for the Hausdorff dimension of products of exact approximation sets. Wang and Wu \cite[Theorem 1.3]{Wang24} established a general principle for estimating the Hausdorff dimension of the Cartesian product of limsup sets generated by balls. As a direct consequence, they obtained the following result.
\begin{theorem}[\cite{Wang24}]
 For each $1 \leq i \leq d$, let $\psi_i:\mathbb{N}\to \mathbb{R}^{+}$ be a positive function. For each $b_i\in\mathbb{N}_{\geq2}$, define
\[
S(\psi_i)= \big\{ x_i \in [0,1] : \parallel b_i^nx_i\parallel<\psi(n), \text{ for infinitely many } n \in \mathbb{N} \big\}.
\]
Then
$$\dim_\mathcal{H} (S(\psi_1) \times \cdots \times S(\psi_d)) \geq \min \big\{ d-1 + \dim_\mathcal{H}S(\psi_i) : 1 \leq i \leq d \big\}.$$
\end{theorem}
Subsequently, Cheng~\cite{Cheng25} extended this result to the path-dependent shrinking target sets in the dyadic and triadic systems.

%Our first result describes the Hausdorff-$f$ measure of $E_{\beta}(\psi)$.
%\begin{theorem}
%Let $\psi : \mathbb{R}_{>0} \to \mathbb{R}_{>0}$ and $f$ be a gauge function.
%Then
%\[
%\mathcal{H}^f(E_{\beta}(\psi)) =
%\begin{cases}
%0, & \text{if } \sum_{n=1}^{\infty} n f(\psi(n)) < \infty, \\[0.5em]
%\mathcal{H}^f([0,1]), & \text{if } \sum_{n=1}^{\infty} n f(\psi(n)) = \infty,
%\end{cases}
%\]
%where $\mathcal{H}^f$ denotes the $f$-dimensional Hausdorff measure.
%\end{theorem}

%Our last result tells us that the Hausdorff dimension of product of $E_{\beta_i}(\psi_i)$.
Motivated by the work of Zhang and Zhong~\cite{Zhang24}, we study the metrical theory of the Cartesian product of exact approximation sets in $\beta$-expansions.
\begin{theorem}\label{mainthm}
Let $d\geq2$ be an integer. For each $1 \leq i \leq d$, let $\beta_i>1$ be real numbers and let $\psi_i : \mathbb{R}_{>0} \to \mathbb{R}_{>0}$
be non-increasing functions with $\psi_i(n)\rightarrow0$ as $n\rightarrow\infty$. Then we have
\begin{align*}
\dim_\mathcal{H} \big( E_{\beta_1}(\psi_1) \times \cdots \times E_{\beta_d}(\psi_d) \big)
&= \min_{1 \leq i \leq d} \Big\{ d - 1 + \dim_\mathcal{H} E_{\beta_i}(\psi_i) \Big\}\\
&=\min_{1 \leq i \leq d} \Big\{ d - 1 + \frac{1}{1+\alpha_i} \Big\}, \quad \alpha_i = \liminf\limits_{n \to \infty} \frac{-\log_{\beta_i} \psi_i(n)}{n}.
\end{align*}
%where $\alpha_i = \liminf\limits_{n \to \infty} \frac{-\log_{\beta_i} \psi_i(n)}{ n}$.
\end{theorem}
%We can easily get the following result.
%\begin{Corollary}
%Let $d\geq2$ be an integer. For each $1 \leq i \leq d$, let $\beta_i>1$ be real numbers and let $\psi_i : \mathbb{R}_{>0} \to %\mathbb{R}_{>0}$
%be non-increasing functions with $\psi_i(n)\rightarrow0$ as $n\rightarrow\infty$. Then we have
%\[
%\dim_\mathcal{H} \big( W_{\beta_i}(\psi_1) \times \cdots \times W_{\beta_d}(\psi_d) \big)
%= \min_{1 \leq i \leq d} \Big\{ d - 1 + \dim_\mathcal{H} W_{\beta_i}(\psi_i) \Big\}
%\]
%where $\lambda_i = \liminf\limits_{n \to \infty} \frac{-\log_{\beta_i} \psi_i(n)}{ n}$.
%\end{Corollary}

\section{Preliminaries}
\subsection{Hausdorff dimension}
In this section, we will provide a definition for the Hausdorff dimension and its properties. We refer the readers to \cite{Falconer14, Mattila95} for further details.

For any set $E \subset \mathbb{R}^d$ and any $\delta > 0$, let $\{U_i\}$ be a countable collection of sets satisfying $|U_i| \leq \delta$ and $E \subset \bigcup\limits_i U_i$.
Let $s \ge 0$ be a real number, and define
\[
\mathcal{H}_{\delta}^s(E)
= \inf\left\{ \sum_i |U_i|^s : \{U_i\} \text{ is a } \delta\text{-cover of } E \right\},
\]
where the infimum is taken over all possible $\delta$-covers of $E$.
The $s$-dimensional Hausdorff measure of $E$ is then defined by
\[
\mathcal{H}^s(E) = \lim_{\delta \to 0} \mathcal{H}_{\delta}^s(E),
\]
and the Hausdorff dimension of $E$ by
\[
\dim_{\mathcal{H}} E
= \inf\{ s \ge 0 : \mathcal{H}^s(E) = 0 \}
= \sup\{ s \ge 0 : \mathcal{H}^s(E) = \infty \}.
\]

%\begin{pro}
%If $E \subset F$, then $\dim_{\mathcal{H}}(E) \leq \dim_{\mathcal{H}}(F)$.
%Furthermore, if $\{E_i\}_{i\ge 1}$ is a countable collection of subsets of $\mathbb{R}$, then
%\[
%\dim_{\mathcal{H}}\!\left(\bigcup_{i=1}^\infty E_i\right)
%= \sup_{i \ge 1} \dim_{\mathcal{H}} E_i.
%\]
%\end{pro}

The following result provides a general method for estimating lower bounds of Hausdorff dimensions, and is commonly known as the \emph{Mass Distribution Principle}.

\begin{proposition}[Mass Distribution Principle \cite{Falconer14}]\label{Mdp}
Let $E$ be a Borel measurable subset of $\mathbb{R}^d$, and let $\mu$ be a Borel measure with $\mu(E) > 0$.
Assume that there exist positive constants $c$ and $\delta$ such that for all $x \in \mathbb{R}^d$ and all $r \in (0,\delta)$,
\[
\mu(B(x,r)) \leq c\, r^s.
\]
Then $\mathcal{H}^s(E) \geq \frac{\mu(E)}{c}$
and hence $\dim_{\mathcal{H}} E \geq s$.
\end{proposition}

\subsection{$\beta$-expansion.}
We begin with a brief overview of some fundamental properties of $\beta$-expansions of real numbers and establish some necessary notation.

%For any $\beta>1$, let $T_{\beta}$ be the $\beta$-transformation defined on $[0,1]$ by $$T_{\beta}(x)=\beta x-\lfloor\beta x\rfloor,$$ where $\lfloor\cdot\rfloor$ denotes the integer part of a real number. Then for every real number $x\in [0,1]$ can be expressed uniquely as a finite or infinite series
%\begin{equation}\label{e1}
%x=\frac{\epsilon_1(x,\beta)}{\beta}+
%\frac{\epsilon_2(x,\beta)}{\beta^2}+\cdots+\frac{\epsilon_n(x,\beta)}{\beta^n}+\frac{T_{\beta}^{n}(x)}{\beta^n}.
%\end{equation}
% we have
%$$x=\sum_{i\geq 1} \frac{\epsilon_i(x, \beta)}{\beta^i}, $$
%where $\epsilon_n(x, \beta)=\lfloor\beta T^{n-1}_\beta x\rfloor$ for $n\geq 1.$
%The series (\ref{e1}) is called the $\beta$-expansion of $x$ and $\{\epsilon_n(x,\beta)\}_{n\geq1}$ the $n$-th digit of $x.$ We also write (\ref{e1}) as $$x=(\epsilon_1(x,\beta),\cdots,\epsilon_n(x,\beta),\cdots).$$
By the definition of $T_\beta$, it is clear that, for $n \ge 1$,
$\varepsilon_n(x, \beta)$ belongs to the alphabet
$A = \{0, 1, \dots, \lceil \beta - 1 \rceil\}$,
where $\lceil y \rceil$ denotes the smallest integer greater than or equal to $x$.
When $\beta$ is not an integer, not every sequences in $A^{\mathbb{N}}$
corresponds to the $\beta$-expansion of some $x \in [0,1]$.
This leads to the notion of $\beta$-admissible sequences.
\begin{definition}   A finite or an infinite sequence $(\varepsilon_1,\cdots, \varepsilon_n, \cdots)$ is called $\beta$-admissible, if there exists an $x\in [0,1]$ such that the $\beta$-expansion of $x$ begins with $\varepsilon_1,\cdots, \varepsilon_n,\cdots$.
\end{definition}
We denote  by $\Sigma_{\beta}^n$ the set of all $\beta$-admissible sequences of length $n$ and
 by  $\Sigma_{\beta}$ the set of all infinite admissible sequences.
% $$\Sigma_{\beta} = \Big\{(\varepsilon_1,\cdots, \varepsilon_n,\cdots)\in \mathcal{A}^{\mathbb{N}}: \text{there exists} \ x\in[0,1] \text{such\ that}\ \varepsilon(x)=(\varepsilon_1,\cdots, \varepsilon_n,\cdots) \Big\},$$}
 %where
%\begin{equation*}
%\mathcal{A}=\left\{
%\begin{aligned}
%& \big\{0, 1, \ldots, \beta-1 \big\} \quad\quad\quad\quad\;\,\,\text{if $\beta$ is an integer}\ ;\\
%& \big\{0, 1, \ldots, \lfloor\beta\rfloor\big\}
% \quad\quad\quad\quad\quad\;\,  \text{if $\beta$ is not an integer}.
%\end{aligned}\right.
%\end{equation*}

The infinite $\beta$-expansion of $1$ is central to the study of $\beta$-admissible sequences.
Denote it by
\[
\varepsilon(1, \beta) := (\varepsilon_1(1, \beta), \varepsilon_2(1, \beta), \dots).
\]
If there are infinitely many $n$ with
$\varepsilon_n \neq 0$, we say that $\varepsilon(1, \beta)$ is infinite.
Otherwise, there exists $m \in \mathbb{N}$ such that
$\varepsilon_m(1, \beta) \neq 0$ and $\varepsilon_n(1, \beta) = 0$ for all $n > m$. Such $\beta$ is called a {\em simple Parry number}, and
$\varepsilon(1, \beta)$ is finite with length $m$.

Define the sequence $\varepsilon^*(1, \beta):=(\epsilon_1^*(\beta), \dots, \varepsilon_n^*(\beta),\dots)$ by
\[
\varepsilon^*(1, \beta) :=
\begin{cases}
\varepsilon(1, \beta), & \text{if}\ \varepsilon(1, \beta)\ \text{is infinite},\\[1mm]
(\varepsilon_1(1, \beta), \dots, \varepsilon_{m-1}(1, \beta), \varepsilon_m(1, \beta)-1)^\infty, &
\text{if}\ \varepsilon(1, \beta)\ \text{is finite with length $m$},
\end{cases}
\]
where $(w)^\infty$ denotes the infinite periodic sequence $(w, w, w, \dots)$.
Thus, $\varepsilon^*(1, \beta)$ is always an infinite sequence, and we have
\[
1 = \sum_{n=1}^{\infty} \frac{\varepsilon_n^*(\beta)}{\beta^n}.
\]

Let $\prec$ and $\preceq$ denote the lexicographic order on $A^{\mathbb{N}}$.
Specifically, $\varepsilon \prec \epsilon'$ means that there exists $k \ge 1$ such that
$\epsilon_j = \varepsilon_j'$ for $1 \le j < k$, while $\varepsilon_k < \varepsilon_k'$.
The notation $\varepsilon \preceq \epsilon'$ means that $\varepsilon \prec \varepsilon'$ or
$\varepsilon = \varepsilon'$. This order extends naturally to finite sequences by
padding with zeros:
\[
(\varepsilon_1, \dots, \varepsilon_n) \mapsto (\varepsilon_1, \dots, \varepsilon_n, 0, 0, \dots).
\]

The following result due to Parry \cite{Parry60} provides a criterion
for determining whether a sequence is $\beta$-admissible.
\begin{lemma}[{\cite{Parry60}}]
Let $\beta > 1$ be a real number.
Then a nonnegative integer sequence $\varepsilon = (\varepsilon_1, \varepsilon_2, \dots)$
is $\beta$-admissible if and only if, for any $k \ge 1$,
\[
(\varepsilon_k, \varepsilon_{k+1}, \dots) \prec (\varepsilon_1^*(\beta), \varepsilon_2^*(\beta), \dots).
\]
\end{lemma}

The number of $\beta$-admissible words of length $n$ satisfies the following bounds due to R\'enyi.
 \begin{lemma}[\cite{Renyi57}]\label{Lemma-2.3}
 Let $ \beta>1 $. For any $ n\ge 1 $, one has
	\[\beta^n\le \# \Sigma_\beta^n\le \frac{\beta^{n+1}}{\beta-1},\]
	where ``$\#(\cdot)$'' denotes the cardinality of a finite set.
 \end{lemma}

For any $(\varepsilon_{1}, \cdots, \varepsilon_{n})\in\Sigma_{\beta}^n$,
we call
$$
I_{n}(\varepsilon_{1}, \cdots, \varepsilon_{n}):=\{x\in [0,1):\varepsilon_{j}(x)=\varepsilon_{j},1\leq j\leq n\}
$$
an $n$-th order cylinder. It is clear that $|I_{n}(\varepsilon_{1}, \cdots, \varepsilon_{n})|\leq \beta^{-n}$, where $ |I| $ denotes the diameter of $I$. It is well known that $([0,1], T_\beta)$ is generally not a sub-shift of finite type, which complicates the metric theory of $\beta$-expansions. The main difficulty is the absence of a uniform lower bound for the length of $n$-level cylinders, which can be much smaller than $\beta^{-n}$. To overcome this, we introduce the following notion.
%Moreover, the $n$-th order cylinder $I_n(\varepsilon_1,\cdots,\varepsilon_n)$ is a left-closed and right-open interval of length at most $\beta^n$ with the left endpoint $\frac{\varepsilon_1}{\beta}+\cdots+\frac{\varepsilon_n}{\beta^n}.$
%Note that the unit interval can naturally be partitioned into a disjoint union of cylinders: for any $n\ge 1$, \begin{equation*}\label{f2}[0,1)=\bigcup\limits_{(\varepsilon_1,\cdots,\varepsilon_n)\in\Sigma_{\beta}^n}I_{n}
%(\varepsilon_1,\cdots,\varepsilon_n).\end{equation*}
% Furthermore, the $n$th order cylinder $I_n(\epsilon_1,\cdots,\epsilon_n)$ is a left-closed and right-open interval of length at most $\beta^n$ with the left endpoint $$
%\frac{\epsilon_1}{\beta}+\cdots+\frac{\epsilon_n}{\beta^n}.$$
%So, if a cylinder $I_{n}(\epsilon_{1}\cdots\epsilon_{n})$ of order $n$ is of length $\beta^{-n}$, it is called a full cylinder.
 %its length verifies: $$I_{n}(\epsilon_{1}\cdots\epsilon_{n})|=\frac{1}{\beta^{n}}$$
\begin{definition}
A  cylinder $ I_{n}(\varepsilon_1,\cdots,\varepsilon_n)$ or a sequence $(\varepsilon_1,\cdots,\varepsilon_n)\in\Sigma_\beta^n $ is called full if it has maximal length, that is,
	\[|I_{n}(\varepsilon_1,\cdots,\varepsilon_n)|=\frac{1}{\beta^n}.\]
\end{definition}
Since full cylinder play an important role in the metric properties of $\beta$-expansions. We now introduce several results about the distribution of full cylinders.
\begin{proposition}[\cite{Fan12}]\label{prop1}
An n-th cylinder $I_n(\varepsilon_1,\cdots,\varepsilon_n)$ is full if and only if for any $\beta$-admissible sequence $(\varepsilon'_1,\cdots,\varepsilon'_m) \in \Sigma_\beta^m$ with $m\geq1$, then $(\varepsilon_1,\cdots,\varepsilon_n, \varepsilon'_1,\cdots,\varepsilon'_m)$ is still $\beta$-admissible. Moreover,
$$|I_{n+m}(\varepsilon_{1}, \varepsilon_{2},\ldots, \varepsilon_{n}, \varepsilon_{1}', \varepsilon_{2}',\ldots, \varepsilon_{m}')|=|I_{n}(\varepsilon_{1}, \varepsilon_{2},\ldots, \varepsilon_{n})|\cdot\big|I_{m}(\varepsilon_{1}', \varepsilon_{2}',\ldots, \varepsilon_{m}')|.$$
\end{proposition}

Let $\beta > 1$ and $\varepsilon^*(1,\beta)$ be the $\beta$-expansion
of $1$. For any $n\geq1,$ we denote by $\ell_{n}(1, \beta)$ the longest length of strings of zeroes behind the $n$-th digit in the $\beta$-expansion of 1, that is,
$$\ell_{n}(1, \beta)=\sup\big\{i\geq0: \varepsilon_{n+1}^{\ast}(1)=\cdots=\varepsilon_{n+i}^{\ast}(1)=0 \big\}.$$
The number $\ell_{n}(1, \beta)$ can be used to describe some admissible words.

For $1\leq j\leq d,$ let $\ell_j(1):=\ell_{1}(1, \beta_j).$ Consequently, by Proposition \ref{prop1}, we have
\begin{equation}\label{eq21}
10^{\ell_j(1)+1}\ \textmd{is\ a\ full\ word.}
\end{equation}
\begin{lemma}[\cite{Bugeaud14}]\label{Lemma-2.7}
For $n\geq1$, among every $n+1$ consecutive cylinders of order n, there exists at least one full cylinder.
\end{lemma}

As a consequence, Wang \cite{Wang18} give a corollary concerning the relationship between balls and cylinders.
%\begin{corollary}[Covering property\cite{Bugeaud14}]
%Let $\beta>1$. For any $y\in [0,1]$ and any positive integer $l$, the ball $B(y,\beta^{-l})$ can be covered by at most $4(l+1)$ cylinders of order $l$.
%\end{corollary}
\begin{corollary}[\cite{Wang18}]\label{p23}
Let $J \subset [0,1]$ be an interval. Then for any integer $n$
with $(n+1)\beta^{-n} < |J|$, there exists a full cylinder $I_{n}$ contained in $J$.
\end{corollary}

\section{Proof of theorem \ref{mainthm}}
The upper bound is trivial. Since $E_{\beta_{j}}(\psi_j)\subset[0,1]$ for every $1\leq j\leq d$, the definition of Cartesian product  together with Theorem \ref{thzhang} yields
\begin{align*}
\dim_\mathcal{H} \big( E_{\beta_1}(\psi_1) \times \cdots \times E_{\beta_d}(\psi_d) \big)
&\leq \min_{1 \leq j \leq d} \Big\{ d - 1 + \dim_\mathcal{H} E_{\beta_j}(\psi_j) \Big\}\\
&=\min_{1 \leq j \leq d} \Big\{ d - 1 + \frac{1}{1+\alpha_j} \Big\}.
\end{align*}

We now turn to the lower bound of the Hausdorff dimension of $E_{\beta_1}(\psi_1) \times \cdots \times E_{\beta_d}(\psi_d)$.
To obtain it, we construct a Cantor subset of this product set. We begin by examining the structure of each individual set $E_{\beta_j}(\psi_j)$.

%It is easy to see that the requirement
%\[
%x - \omega_n(x) < \frac{\psi_j(n)}{\beta_j^n}
%\quad \text{for infinitely many } n \in \mathbb{N}
%\]
%is not difficult to fulfill.
%The real challenge is to ensure that for any $c < 1$,
%\[
%x - \omega_n(x) > c\cdot \frac{\psi_j(n)}{\beta_j^n}
%\quad \text{ultimately}.
%\]

Achieving the exact order of approximation requires careful control of the digits, in particular avoiding long zero blocks. We must also verify that the constructed points indeed lie in $E_{\beta_j}(\psi_j)$.

We introduce some notation.
\begin{equation*}
    \begin{split}
   %\bullet &~\text{For}~1\leq j\leq d, l_{1,j}=l_{1,j}(1, \beta_{j})~\text{denote the longest length of strings of zeroes behind the first}\\
%   &\text{ digit in the}~\beta_{j}\text{-expansion of 1};\\
   \bullet &~ a\asymp b~\text{if}~ c^{-1}<a/b<c,~\text{and if}~
a\lesssim b~\text{if}~a\leq cb ~\text{for some constant}\ c>1;\\
\bullet &~ \text{Throughout, numbers carry a double subscript, with the second indicating the}\\
&\text{corresponding }\beta, \text{ while words are marked by a superscript. For instance, }
\delta_{i,j}>0\\
&\text{corresponds to }\beta_j,\ \text{and } \tau_q^{(j)}\in\Sigma_{\beta_j}
\text{ is a word in the }\beta_j\text{-shift space}.
   \end{split}
\end{equation*}

Let $s := \min\limits_{1\leq j\leq d}\Bigl(d-1 + \tfrac{1}{1+\alpha_j}\Bigr)$.
For $n \geq 1$, define
\[
\mathcal{F}_{n, j} = \left\{ \xi \in \Sigma^n_{\beta_j} : I_n(\xi) \ \text{is a full cylinder} \right\}
\]
and
\[
\mathcal{R}_{n, j} = \left\{ \xi \in \Sigma^n_{\beta_j} : I_n(\xi) \ \text{is a full cylinder starting with a non-zero digit} \right\}.
\]
By \eqref{eq21},
\begin{equation}\label{Eq-1}
\mathcal{R}_{n,j} \supset \{\, 10^{\,\ell_j(1)+1}\xi : \xi \in \mathcal{F}_{n-(\ell_j(1)+2),j} \,\}.
\end{equation}

Fix $0 < \eta < 1$, choose $M_0$ sufficiently large such that for any $M \geq M_0$,
\begin{equation}\label{Eq-2}
\frac{\beta_j^{M-\ell_j(1)-2}}{{M-\ell_j(1)-1}} \geq \beta_j^M (1-\eta),
\quad 1\leq j \leq d.
\end{equation}
Combining  Lemma \ref{Lemma-2.3}, Lemma \ref{Lemma-2.7} together with formulas \eqref{Eq-1} and \eqref{Eq-2}, we have
\begin{equation}\label{Eq-M}
\# \mathcal{R}_{M,j} \geq \beta_j^M (1-\eta).
\end{equation}
We now fix an integer $M$ for which \eqref{Eq-M} holds.

To proceed with the construction, we need the following Claim.

\noindent{\bf Claim.}
For each $1\le j\le d$, fix a sequence $\{\delta_{i,j}\}_{i\ge 1}\subset\mathbb{R}_{>0}$ with $\delta_{i,j}\to0$ as $i\to\infty$.
Since $\psi_j(n)\to0$, the definition of $\alpha_j$ implies the existence of a sparse subsequence $\{n'_{i,j}\}_{i\ge 1}$ of $\mathbb{N}$ such that
\begin{equation}\label{Eq-3}
\lim\limits_{i\to\infty}\frac{-\log_{\beta_j}\psi_j(n'_{i,j})}{n'_{i,j}}=
\liminf\limits_{i\to\infty}\frac{-\log_{\beta_j}\psi_j(n)}{n}
= \alpha_j
\end{equation}
and
\begin{equation}\label{Eq-4}
\frac{\log_{\beta_j}\psi_j(n'_{i,j})}{\log\delta_{i,j}}\to\infty,~\text{as}~i\to\infty,
\end{equation}
for example by taking
$\delta_{i,j}\ge\psi_j(n'_{i,j})^{1/i}.$

Recall that $\ell_j(1)$ is the longest length of strings of zeroes behind the first digit in the $\beta_j$-expansion of 1, we can ensure that
\begin{equation}\label{Eq-4.5}
\psi_j(n'_{i,j})
<\beta_j^{-\ell_j(1)-3},
\qquad
\psi_j(n'_{i,j})<\beta_j^{-2M}.
\end{equation}

Let $k_{i,j}\in\mathbb N$ be the integer such that
\begin{equation}\label{Eq-5}
(k_{i,j}+1)\cdot\beta_j^{-k_{i,j}}
\le
\delta_{i,j}\cdot\psi_j(n'_{i,j})
\le
k_{i,j}\cdot\beta_j^{-k_{i,j}+1}.
\end{equation}
By Corollary~\ref{p23}, there exists a full cylinder
\[
I_{k_{i,j}}(\sigma_i^{(j)})
\subset
\bigl((1-\delta_{i,j})\psi_j(n'_{i,j}),\ \psi_j(n'_{i,j})\bigr).
\]
In particularly, \eqref{Eq-4} implies that
\begin{equation}\label{Eq-5.5}
k_{i,j}
  \asymp \log_{\beta_j}\frac1{\delta_{i,j}\psi_j(n'_{i,j})}
  \asymp \log_{\beta_j}\frac1{\psi_j(n'_{i,j})}.
\end{equation}

Let $t_{i,j}$ be the integer satisfying
\begin{equation}\label{Eq-6}
\beta_j^{-(t_{i,j}+1)}
\le
\psi_j(n'_{i,j})
< \beta_j^{-t_{i,j}}.
\end{equation}
Since every element of $I_{k_{i,j}}(\sigma_i^{(j)})$ is smaller than $\psi_j(n'_{i,j})$, the word $\sigma_i^{(j)}$ begins with at least $t_{i,j}$ zeros. By \eqref{Eq-4.5}, $t_{i,j}\ge\ell_j(1)+3$, hence
\begin{equation}\label{Eq-z}
\sigma_i^{(j)}=0^{\,\ell_j(1)+1} z_i^{(j)},
\end{equation}
where $z_i^{(j)}$ is a full word.

%Finally, \eqref{Eq-5-compact} implies
%\[
%(k_{q,h}+1)\beta_h^{-k_{q,h}}
%<\beta_h^{-t_{q,h}},
%\qquad \text{so} \qquad
%k_{q,h}>t_{q,h}+\log_{\beta_h}(k_{q,h}+1)>t_{q,h}.
%\]

\subsection{Construction of the Cantor set.}
We construct the Cantor subset alternately on each component.

%\textcolor{red}{For $1\leq i\leq d$, we choose $p_{i,i}\in\mathbb{N}$  and $1\leq r_{i,i} <M$ such that $n_{i,i} = p_{i,i}M + r_{i,i}$.}

{\em The first level of the Cantor set.}
We choose $p_{1,1}\in\mathbb{N}$  and $1\leq r_{1,1} <M$ such that $n_{1,1}:=n'_{1,1} = p_{1,1}M + r_{1,1}$.
For every $\bos\xi^{(1)}_1=(\xi^{(1)}_1, \ldots, \xi^{(1)}_{p_{1,1}})\in \mathcal{R}^{p_{1,1}}_{M,1}$,
we obtain a full cylinder of the form
\[
I_{n_{1,1}-1}(\bos\xi^{(1)}_{1},0^{r_{1,1}-1}).
\]

Define
\[
E_{1,1} := \Bigl\{ x \in [0,1) : (1-\delta_{1,1})\psi_1(n_{1,1}) < \beta_1^{n_{1,1}}(x - \omega_{n_{1,1}}(x)) < \psi_1(n_{1,1}) \Bigr\}.
\]
By the {\bf Claim}, there exists a full word of the form \eqref{Eq-z} satisfying
\[
I_{k_{1,1}}(\sigma^{(1)}_1) \subset \bigl((1-\delta^{(1)}_1)\psi_1(n_{1,1}),\, \psi_1(n_{1,1})\bigr).
\]
This gives the full cylinder
\[
I_{n_{1,1}+k_{1,1}}(\bos\xi^{(1)}_{1}, 0^{\,r_{1,1}-1}, 1, \sigma^{(1)}_1)
\subset E_{1,1}.
\]

%We claim that the cylinder $I_{n_{1,1}+k_{1,1}}(\bos\xi^{(1)}_{1}, 0^{\,r_{1,1}-1}, 1, \sigma^{(1)}_1)$
%is a full cylinder. By Lemma~2.7 (2), it suffices to show that
%$I_{1+k_{1,1}}(1,\sigma^{(1)}_1)$ is a full cylinder.
%Recalling \eqref{Eq-z}, one has
%\[
%(1,\sigma^{(1)}_1) = 10^{\ell_{1}(1)+1}z^{(1)}_1.
%\]
%Since both $10^{\ell_1(1)+1}$ and $z^{(1)}_1$ are full words, it follows that $(1,\sigma^{(1)}_1)$
%is also a full word. Thus, by Proposition \ref{prop1}, the cylinder $I_{n_{1,1}+k_{1,1}}(\bos\xi^{(1)}_{1}, 0^{\,r_{1,1}-1}, 1, \sigma^{(1)}_1)$
%is a full cylinder.

Set $N_{1,1}:=n_{1,1}+k_{1,1}$ and define
\[
\mathcal{C}_{1} := \Bigl\{ I_{N_{1,1}}(\bos\xi^{(1)}_{1}, 0^{\,r_1-1}, 1, \sigma^{(1)}_{1})\times[0,1)^{d-1}
: \bos\xi^{(1)}_{1} \in \mathcal{R}^{p_{1,1}}_{M,1} \Bigr\},
\]
then $\# \mathcal{C}_1=(\# \mathcal{R}_{M,1})^{p_{1,1}}.$

Now we cut these rectangles into approximate squares. We write $\Gamma^{(1)}_{1}$ for any word of the form $(\bos\xi^{(1)}_{1}, 0^{\,r_1-1}, 1, \sigma^{(1)}_{1})$ in $\mathcal{C}_{1}.$ Let $n_{1,j}=\lfloor(n_{1,1}+k_{1,1})\log_{\beta_j}\beta_1\rfloor$ and %For any $I_{N_{1,1}}(\Gamma^{(1)}_{1})\times[0,1)^{d-1}\in \mathcal{C}_{1}$,
writing $N_{1,j}:=n_{1,j}=p_{1,j}M+r_{1,j}$ with $1\le r_{1,j}<M$, we define
\[
\mathcal{B}_1:=
\Bigl\{
I_{N_{1,1}}(\Gamma^{(1)}_{1})\times
\prod_{j=2}^d I_{N_{1,j}}(\boldsymbol{\xi}^{(j)}_{1},0^{r_{1,j}}):
\boldsymbol{\xi}^{(j)}_{1}\in\mathcal{R}^{p_{1,j}}_{M,j}
\Bigr\},
\]
and thus
$$\#\mathcal{B}_1=\prod_{j=1}^d (\#\mathcal{R}_{M,j})^{p_{1,j}}.$$
It should be noted that the word $\Gamma^{(1)}_{1}$ and the words $(\boldsymbol{\xi}^{(j)}_{1},0^{r_{1,j}})$ $(2\leq j\leq d)$ have different forms, but for convenience, we still denote $(\boldsymbol{\xi}^{(j)}_{1},0^{r_{1,j}})$ by $\Gamma^{(j)}_{1}$ for $2\leq j\leq d$. The same convention will be used for the subsequent levels as well.

%And we emphasize the distinction between $\Gamma_1^{(1)}$ and $\Gamma_1^{(j)}$ ($1<j\leq d$).}

{\em The second level of the Cantor set.}
Fix an element $J_{1}:=\prod\limits_{1\leq j\leq d}I_{N_{1,j}}(\Gamma^{(j)}_{1})\in \mathcal{B}_{1}$.
We choose a large integer $p_{2,2}\in\mathbb{N}$ and $1\leq r_{2,2} <M$ such that
$$\frac{N_{1,2}}{p_{2,2}M+r_{2,2}}\leq \eta$$
and put $n_{2,2} = p_{2,2}M + r_{2,2}$.
For every $\bos\xi^{(2)}_2=(\xi^{(2)}_1, \ldots, \xi^{(2)}_{p_{2,2}})\in R^{p_{2,2}}_{M,2}$.
We obtain a full cylinder of the form
\[
I_{N_{1,2}+n_{2,2}-1}(\Gamma^{(2)}_{1},\bos\xi^{(2)}_{2},0^{\,r_{2,2}-1}).
\]

Define
\[
E_{2,2} := \Bigl\{ x \in [0,1) : (1-\delta_{2,2})\cdot\psi_2(N_{1,2}+n_{2,2}) < \beta_2^{N_{1,2}+n_{2,2}}(x - \omega_{N_{1,2}+n_{2,2}}(x)) < \psi_2(N_{1,2}+n_{2,2}) \Bigr\}.
\]
By the {\bf Claim}, there exists a full word of the form \eqref{Eq-z} satisfying
\[
I_{k_{2,2}}(\sigma^{(2)}_{2}) \subset \bigl((1-\delta_{2,2})\psi_2(N_{1,2}+n_{2,2}),\, \psi_2(N_{1,2}+n_{2,2})\bigr).
\]
This gives the full cylinder
\[
I_{N_{1,2}+n_{2,2}+k_{2,2}}(\Gamma^{(2)}_{1},\bos\xi^{(2)}_{2},0^{\,r_{2,2}-1}, 1, \sigma^{(2)}_2)
\subset E_{2,2}
\]

%We claim that the cylinder $I_{N_{1,2}+n_{2,2}+k_{2,2}}(\Gamma^{(2)}_{1},\bos\xi^{(2)}_{2},0^{\,r_{2,2}-1}, 1, \sigma^{(2)}_2)$
%is a full cylinder. By Lemma~2.7 (2), it suffices to show that
%$I_{1+k_{2,2}}(1,\sigma^{(2)}_2)$ is a full cylinder. Recalling (13), one has
%\[
%(1,\sigma^{(2)}_2) = 10^{\ell_2(1)+1}z^{(2)}_{2}.
%\]
%Since both $10^{\ell_2(1)+1}$ and $z^{(2)}_{2}$ are full words, it follows that $(1,\sigma^{(2)}_2)$
%is also a full word. This proves the claim.

Set $N_{2,2}:=N_{1,2}+n_{2,2}+k_{2,2}$, we define
\[
\mathcal{C}_{2}(J_1) := \left\{ I_{N_{2,2}}(\Gamma^{(2)}_{1},\bos\xi^{(2)}_{2},0^{\,r_{2,2}-1}, 1, \sigma^{(2)}_2)\times\prod_{j\ne 2}I_{N_{1,j}}(\Gamma^{(j)}_{1})
: \bos\xi^{(2)}_{2} \in R^{p_{2,2}}_{M,2} \right\},
\]
then
$$\# \mathcal{C}_2(J_1)=\#\mathcal{B}_1\cdot(\# R_{M,2})^{p_{2,2}}.$$
We abbreviate by $\Gamma^{(2)}_{2}$ a generic word $(\Gamma^{(2)}_{1},\bos\xi^{(2)}_{2},0^{\,r_{2,2}-1}, 1, \sigma^{(2)}_2)$.

Now we cut these rectangles into approximate squares. For any $j\ne 2$, take $n_{2,j}=\lfloor(n_{2,2}+k_{2,2})\log_{\beta_j}\beta_2\rfloor$ and write $n_{2,j}=p_{2,j}M+r_{2,j}$ with
 $1\le r_{2,j}<M$.
%For any $I_{N_{2,2}}(\Gamma^{(2)}_{2})\times\prod_{j\ne 2}I_{N_{1,j}}(\Gamma^{(j)}_{1})\in \mathcal{C}_{2}$, s
Let $N_{2,j}:=N_{1,j}+n_{2,j}$ and define
$$\mathcal{B}_{2}(J_1):=\left\{I_{N_{2,2}}(\Gamma^{(2)}_{2})\times \prod_{ j\ne2}I_{N_{2,j}}(\Gamma_{1}^{(j)},\bos\xi^{(j)}_{2},0^{r_{2,j}}): \bos\xi^{(j)}_{2}\in R^{p_{2,j}}_{M,j} \right\}.$$
One has
$$\# \mathcal{B}_{2}(J_1)=\#\mathcal{B}_1\cdot\prod_{1\leq j \leq d} (\# R_{M,j})^{p_{2,j}}.$$

Therefore, the second level of the Cantor set is defined as
$$\mathcal{C}_{2}=\bigcup_{J\in \mathcal{B}_1}\mathcal{C}_2(J),\quad \mathcal{B}_{2}=\bigcup_{J\in \mathcal{B}_1}\mathcal{B}_2(J).$$

%\textcolor{red}{As before, for convenience, we still denote $(\Gamma_{1}^{(j)},\bos\xi^{(j)}_{2},0^{r_{2,j}})$ by $\Gamma_2^{(j)}$ for $j\ne2$. }

{\em The $q$ level of the Cantor set.}
Let $q:=ld+h$ for some integer $l\geq0$ and $1\leq h \leq d-1$.
Suppose that the collection $\mathcal{B}_{q-1}$ has already been defined, and let
\[
J_{q-1}:=\prod_{1\leq j \leq d}I_{N_{q-1,j}}(\Gamma^{(j)}_{q-1})\in \mathcal{B}_{q-1}.
\]
Choose $p_{q,h}\in\mathbb{N}$ and $1\leq r_{q,h} <M$ such that
\begin{equation}\label{Eq-eta}
\frac{N_{q-1,h}}{p_{h,h}M+r_{q,h}}\leq \eta,
\end{equation}
and write $n_{q,h} := p_{q,h}M + r_{q,h}$.
For each $\bos\xi^{(h)}_q=(\xi^{(h)}_1, \ldots, \xi^{(h)}_{p_{q,h}})\in \mathcal{R}^{p_{q,h}}_{M,h}$,
we obtain a full cylinder of the form
\[
I_{N_{q-1,h-1}+n_{q,h}-1}\bigl(\Gamma^{(h-1)}_{q-1},\bos\xi^{(h)}_{q},0^{\,r_{q,h}-1}\bigr).
\]

Define
$$
\begin{aligned}
E_{q,h} := \Bigl\{ x \in [0,1) : (1-\delta_{q,h})\psi_h(N_{q-1,h}+n_{q,h})
&< \beta_h^{N_{q-1,h}+n_{q,h}}\bigl(x - \omega_{N_{q-1,h}+n_{q,h}}(x)\bigr)\\
&\qquad\qquad\qquad\qquad< \psi_h(N_{q-1,h}+n_{q,h}) \Bigr\}.
\end{aligned}
$$
By the {\bf Claim}, there exists a full word $\sigma^{(h)}_{q}$ of form \eqref{Eq-z} satisfying
\[
I_{k_{q,h}}(\sigma^{(h)}_{q}) \subset
\bigl((1-\delta_{q,h})\psi_h(N_{q-1,h}+n_{q,h}),\, \psi_h(N_{q-1,h}+n_{q,h})\bigr).
\]
This gives the full cylinder
\[
I_{N_{q-1,h}+n_{q,h}+k_{q,h}}
\bigl(\Gamma^{(h)}_{q-1},\bos\xi^{(h)}_{q},0^{\,r_{q,h}-1}, 1, \sigma^{(h)}_q\bigr)
\subset E_{q,h}.
\]

%We claim that the cylinder
%Like before, the cylinder
%\[
%I_{N_{q-1,h}+n_{q,h}+k_{q,h}}
%(\Gamma^{(h)}_{q-1},\bos\xi^{(h)}_{q},0^{\,r_{q,h}-1}, 1, \sigma^{(h)}_q)
%\]
%is indeed full.
%By Lemma~2.7 (2), it suffices to check that
%$I_{1+k_{q,h}}(1,\sigma^{(q)}_h)$ is full. Recalling (13), we have
%\[
%(1,\sigma^{(h)}_q) = 10^{\ell_{q-1,h}+1}\epsilon^{(h)}_{q}.
%%\]
%Since both $10^{\ell_{q-1,h}+1}$ and $\epsilon^{(h)}_{q}$ are full words,
%it follows that $(1,\sigma^{(q)}_h)$ is also a full word. This proves the claim.

Set $N_{q,h}:=N_{q-1,h}+n_{q,h}+k_{q,h}$, define
\[
\mathcal{C}_{q}(J_{q-1}) := \left\{
I_{N_{q,h}}(\Gamma^{(h)}_{q-1},\bos\xi^{(h)}_{q},0^{\,r_{q,h}-1}, 1, \sigma^{(h)}_q)
\times\prod_{j\ne h}I_{N_{q-1,j}}(\Gamma^{(j)}_{q-1})
: \bos\xi^{(h)}_{q} \in R^{p_{q,h}}_{M,h} \right\}.
\]
Thus
\[
\# \mathcal{C}_q(J_{q-1})=\#\mathcal{B}_{q-1}\cdot(\# R_{M,h})^{p_{q,h}}.
\]
We abbreviate by $\Gamma^{(h)}_{q}$ a generic word $(\Gamma^{(h)}_{q-1},\bos\xi^{(h)}_{q},0^{\,r_{q,h}-1}, 1, \sigma^{(h)}_q).$

Now we cut these rectangles into approximate squares. For each $j\ne h$, set
$n_{q,j}=\Bigl\lfloor (n_{q,h}+k_{q,h})\log_{\beta_j}\beta_h \Bigr\rfloor$,
and write $n_{q,j}=p_{q,j}M+r_{q,j}$ with $r_{q,j}<M$.
%For any element
%\[
%I_{N_{q,h}}(\Gamma^{(h)}_{q})\times\prod_{j\ne h}I_{N_{q-1,j}}(\Gamma^{(j)}_{q-1})\in \mathcal{C}_{q}(J_{q-1}),
%\]
Let $N_{q,j}=N_{q-1,j}+n_{q,j}$ and define
\[
\mathcal{B}_{q}(J_{q-1}):=\left\{
I_{N_{q,h}}(\Gamma^{(h)}_{q})
\times \prod_{j\ne h}I_{N_{q,j}}(\Gamma_{q-1,j},\bos\xi^{(j)}_{q},0^{r_{q,j}})
: \bos\xi^{(j)}_{q}\in R^{p_{q,j}}_{M,j} \right\}.
\]
Then
\[
\# \mathcal{B}_{q}(J_{q-1})=\#\mathcal{B}_{q-1}\prod_{1\leq j \leq d} (\# R_{M,j})^{p_{q,j}}.
\]

Therefore, the $q$-th level of the Cantor set is given by
\[
\mathcal{C}_{q}=\bigcup_{J\in \mathcal{B}_{q-1}}\mathcal{C}_q(J),\quad
\mathcal{B}_{q}=\bigcup_{J\in \mathcal{B}_{q-1}}\mathcal{B}_q(J).
\]

\noindent{\bf The Cantor set.} Continuing this construction, we obtain a nested family $\{\mathcal{B}_q\}_{q\ge1}$, and define the desired Cantor set by
$$\mathcal{B}_{\infty}= \bigcap^{\infty}_{q=1}\bigcup_{J\in \mathcal{B}_{q}}J.$$
\begin{remark}\label{3.1}
For any approximative square in $\mathcal{B}_q$ with $q=ld+h$, we have the following relations:
\begin{itemize}
\item When $j=h$:
  $$
  N_{q,j} = N_{q-1} + n_{q,j} + k_{q,j};
  $$
\item When $j\ne h$:
  $$
  N_{q,j} = N_{q-1} + n_{q,j}.
  $$
\end{itemize}
Since for every $j\neq h$, we have $\beta_j^{n_{q,j}}\asymp \beta_h^{n_{q,h}+k_{q,h}}$. Then
$$
\beta_1^{N_{q,1}} \asymp \beta_2^{N_{q,2}}\asymp\cdots \asymp \beta_d^{N_{q,d}},
$$
\end{remark}
\begin{remark}
Let $q = ld + h$. Although for $1 \leq j \leq d$, we have
\[
n_{q,j} = p_{q,j} M + r_{q,j}, \qquad 1 \leq r_{q,j} \leq M-1.
\]
It should be emphasized that the integer $n_{q,j}$ ($j\neq h$) depends on
$n_{q,h}+k_{q,h}$ according to the relation stated in the Remark~\ref{3.1}. Moreover,
$\delta_{q,j}$, $k_{q,j}$ and $t_{q,j}$ appear only in the case $j=h$.
In addition, it is clear that the sequence $\{N_{i-1,j}+n_{i,j}\}$ defined in the construction still  satisfies the conditions in the {\bf Claim}, hence we indeed deal with the sequence $\{N_{i-1,j}+n_{i,j}\}$.
\end{remark}

We now verify that $\mathcal{B}_{\infty}$ is contained in $\prod^{d}\limits_{j=1}E_{\beta_j}(\psi_j).$
\begin{proposition}
One has
$$\mathcal{B}_{\infty}\subset \prod^{d}_{j=1}E_{\beta_j}(\psi_j).$$
\begin{proof}
Let $\bos x=(x_1,\cdots,x_d) \in \mathcal{B}_{\infty}$. By the definition of the Cartesian product, it suffices to show that $x_j\in E_{\beta_j}(\psi_j)$ for each $1 \le j \le d$.

Fix $0\leq j \leq d-1$.  We consider the quantity
\[
\beta_j^n \big(x_j - \omega_n(x_j)\big).
\]
Let $q = ld + j$ be the integer such that
$N_{q-1,j} + n_{q,j} \le n < N_{q+d-1,j} + n_{q+d,j}$.

\noindent{(I)} If $n = N_{q-1,j}+n_{q,j}$,  by the definition of $E_{q,j}$, we have
\begin{equation}\label{Eq-8}
(1 - \delta_{q,j})\cdot \psi_j(n) < \beta_j^n (x_j - \omega_n(x_j)) < \psi_j(n).
\end{equation}

\noindent{(II)}If $N_{q-1,j}+n_{q,j} < n \leq N_{q,j}=N_{q-1,j}+n_{q,j}+k_{q,j}$, then $\sigma^{(j)}_q$ begins with at least $t_{q,j}$ zeros by \eqref{Eq-z}, then
\begin{equation}\label{Eq-0}
\varepsilon_{N_{q-1,j}+n_{q,j}+1}(x_j) = \cdots = \varepsilon_{N_{q-1,j}+n_{q,j}+t_{q,j}}(x_j) = 0.
\end{equation}
Since $k_{q,j} > t_{q,j} + \log_{\beta_j}(k_{q,j} + 1) > t_{q,j}$,
we divide (II) into two subcases:

\noindent{(i)} If $N_{q-1,j}+n_{q,j} < n \leq N_{q-1,j}+n_{q,j}+t_{q,j}$, then
\begin{align*}
\beta_j^n (x_j - \omega_n(x_j))
&= \beta_j^n (x_j - \omega_{\,N_{q-1,j}+n_{q,j}}(x_j))\quad \text{by  \eqref{Eq-0}}\\
&> \beta_j^{\,n - (N_{q-1,j}+n_{q,j})} \cdot \bigl(1 - \delta_{q,j}\cdot \psi_j(\,N_{q-1,j}+n_{q,j})\bigr) \quad \text{by  \eqref{Eq-8}}\\
&> (1 - \delta_{q,j})\cdot \psi_j(\,N_{q-1,j}+n_{q,j})\\
&> (1 - \delta_{q,j})\cdot \psi_j(n).
\end{align*}

\noindent{(ii)} If $N_{q-1,j}+n_{q,j}+t_{q,j} < n \leq N_{q,j}=N_{q-1,j}+n_{q,j}+k_{q,j}$, since each word in $R_{M,j}$ begins with a non-zero digit, so $\varepsilon_{N_{q,j}+1}(x_j) \neq 0$. Hence,
\begin{equation}\label{Eq-9}
\begin{aligned}
\beta_j^n (x_j - \omega_n(x_j))
&=\beta_j^n\big( \frac{ \varepsilon_{n+1}(x_j)}{\beta_j^{n+1}}
   + \cdots + \frac{\varepsilon_{N_{q,j}+1}(x_j)}{\beta_j^{N_{q,j}+1}} + \cdots\big)\\[6pt]
&\geq \beta_j^n \cdot \frac{1}{\beta_j^{N_{q,j}+1}}\geq \frac{1}{\beta_j^{k_{q,j} - t_{q,j} + 1}}.
\end{aligned}
\end{equation}
By \eqref{Eq-5} and \eqref{Eq-6}, we also have
\[
\delta_{q,j} \cdot \frac{\psi_{j}(N_{q-1,j}+n_{q,j})}{k_{q,j}} \leq \frac{1}{\beta_j^{k_{q,j} - 1}}
\quad \text{and} \quad
\frac{1}{\beta_j^{t_{q,j} - 2}} = \frac{\beta_j^3}{\beta_j^{t_{q,j} + 1}} < \psi_j(N_{q-1,j}+n_{q,j}) \cdot \beta_j^3,
\]
which implies that
\begin{equation}\label{Eq-10}
\frac{1}{\beta_j^{k_{q,j} - t_{q,j} + 1}}
= \frac{1/\beta_j^{k_{q,j} - 1}}{1/\beta_j^{t_{q,j} - 2}}
\geq \beta_j^{-3} \cdot \frac{\delta_{q,j}}{k_{q,j}}
 \overset{\text{by  \eqref{Eq-4.5}}}{\geq} \psi(N_{q-1,j}+n_{q,j}).
\end{equation}
%where the last inequality is from (11).
Combining \eqref{Eq-9} and \eqref{Eq-10}, we have
\[
\beta_j^n (x_j - \omega_n(x_j)) > (1 - \delta_{q,j})\cdot \psi(N_{q-1,j}+n_{q,j}) > (1 - \delta_{q,j})\cdot \psi(n).
\]

\noindent{(III)} If $N_{q,j}\leq n <N_{q+d-1,j}$, suppose $N_{q+i,j}\leq n< N_{q+i+1,j}$ for any $0\leq i <d-2$.
We consider the following two subcases.

\noindent{(i)} If there exists $1 \leq r \leq p_{q+i+1,j}$ such that
\[
N_{q+i,j} + (r - 1)M < n \leq N_{q+i,j} + rM,
\]
since $\varepsilon_{N_{q+i,j} + rM+1}(x_j) \neq 0$ for any $1 \leq r \leq p_{k+1}$, we have
\begin{equation*}
\begin{aligned}
\beta_j^n (x_j - \omega_n(x_j))
&= \beta_j^n\big(\frac{ \varepsilon_{n+1}(x_j)}{\beta_j^{n+1}}
+ \cdots + \frac{\varepsilon_{N_{q+i,j} + rM+1}(x_j)}{\beta_j^{N_{q+i,j} + rM+1}} + \cdots\big)\\
&\geq \beta_j^n \cdot \frac{1}{\beta_j^{N_{q+i,j} + rM+1}}\geq \frac{1}{\beta_j^{M+1}}\\
&> (1 - \delta_{q,j})\cdot \psi(N_{q-1,j}+n_{q,j}) > (1 - \delta_{q,j})\cdot \psi(n).
\end{aligned}
\end{equation*}

\noindent{(ii)} If $N_{q+i,j} + p_{q+i+1,j}M < n < N_{q+i+1,j}=N_{q+i,j}+n_{q+i+1,j}$,
then $\varepsilon_{N_{q+i+1,j}+1}(x_j) \neq 0$ by the definition of $R_{M,j}$ and hence
\begin{equation*}
\begin{aligned}
\beta_j^n (x_j - \omega_n(x_j))
&= \beta_j^n \big(\frac{\varepsilon_{n+1}(x_j)}{\beta_j^{n+1}}
+ \cdots + \frac{\varepsilon_{N_{q+i+1,j}+1}(x_j)}{\beta_j^{N_{q+i+1,j}+1}} + \cdots\big)\\
&\geq \beta_j^n \cdot \frac{1}{\beta_j^{N_{q+i+1,j}+1}}
\geq \frac{1}{\beta_j^{r_{q+i+1,j}+1}}
\geq \frac{1}{\beta_j^{M}}\\
&>( 1 - \delta_{q,j})\cdot \psi(N_{q-1,j}+n_{q,j}) > (1 - \delta_{q,j})\cdot \psi_j(n).
\end{aligned}
\end{equation*}

\noindent{(IV)} If $N_{q+d-1,j}\leq n <N_{q+d-1,j}+n_{q+d,j}$, a similar argument to that in case (III) gives
\[
\beta_j^n (x_j - \omega_n(x_j)) > (1 - \delta_{q,j})\cdot\psi_j(n).
\]
The case (I) shows that there exist at least $l=\bigl\lfloor\tfrac{q}{d}\bigr\rfloor$ solutions of
\[
\beta_j^n(x_j-\omega_n(x_j))<\psi_j(n).
\]
On the other hand, (I)-(IV) show that
\[
\beta_j^n(x_j-\omega_n(x_j))<(1-\delta_{q,j})\psi_j(n)
\]
has no solution for $N_{q-1,j}+n_{q,j}\le n< N_{q+d-1,j}+n_{q+d,j}$. Thus $x_j\in E_{\beta_j}(\psi_j)$.

Repeating this for each $1\le j\le d$ yields
\[
\mathcal{B}_\infty \subset \prod_{j=1}^{d}E_{\beta_j}(\psi_j).
\]
\end{proof}
\end{proposition}

\subsection{Supporting measure}

We now distribute a probability measure $\mu$ on $\mathcal{B}_\infty$.
First, define $\mu([0,1)^{d}) = 1$.
\begin{itemize}[leftmargin=15pt]
    \item {Measure of elements in $\mathcal{B}_1$}. For a generic rectangle $I_{N_{1,1}}(\Gamma^{(1)}_{1}) \times [0, 1)^{d-1} \in \mathcal{C}_1$, define
    \[
    \mu(I_{N_{1,1}}(\Gamma^{(1)}_{1}) \times [0, 1)^{d-1}) = \frac{1}{\# R_{M,1}^{p_{1,1}}},
    \]
The measure of an approximative square $J_{1}=\prod\limits_{1\leq j\leq d}I_{N_{1,j}}(\Gamma^{(j)}_{1})\in\mathcal{B}_1$ is defined as
    \[
\begin{aligned}
    \mu(J_{1})
    &= \prod_{j\ne 1}\frac{1}{\# R_{M,j}^{p_{1,j}}}\cdot \mu(I_{N_{1,1}}(\Gamma^{(1)}_{1}) \times [0, 1)^{d-1}) \\
    &= \prod^{d}_{j=1}\frac{1}{\# R_{M,j}^{p_{1,j}}}.
\end{aligned}
\]

\item {Measure of elements in $\mathcal{B}_q$}. Assume that the measure of approximation squares $J_{q-1} \in \mathcal{B}_{q-1}$ has been well defined. %We define the measure on $\mathcal{B}_q$ by defining the measures on the offsprings of $J_{q-1}$.
Write $q = ld + h$ for some $1 \leq h \leq d-1$. For a generic rectangle
\[
I_{N_{q,h}}(\Gamma^{(h)}_{q}) \times \prod_{j \ne h} I_{N_{q-1,j}}(\Gamma^{(j)}_{q-1}) \in \mathcal{C}_q,
\]
we define
\begin{align*}
\mu\Bigl(I_{N_{q,h}}(\Gamma^{(h)}_{q}) \times \prod_{j \ne h} I_{N_{q-1,j}}(\Gamma^{(j)}_{q-1})\Bigr)
&:= \mu(J_{q-1}) \frac{1}{\# R_{M,h}^{p_{q,h}}}.
\end{align*}
The measure of the approximative square
$
J_q = \prod\limits_{1 \le j \le d} I_{N_{q,j}}(\Gamma^{(j)}_{q}) \in \mathcal{B}_q
$
is defined as
\begin{align*}
\mu(J_q)
&= \mu(J_{q-1})\frac{1}{\# R_{M,h}^{p_{q,h}}}\cdot \prod_{j \ne h} \frac{1}{\# R_{M,j}^{p_{q,j}}} \notag\\
&= \mu(J_{q-1}) \prod_{j=1}^{d} \frac{1}{\# R_{M,j}^{p_{q,j}}}.
\end{align*}
\end{itemize}

It is clear that $\mu$ satisfies the consistency property and can then be uniquely extended to a non-atomic Borel measure supported on $\mathcal{B}_\infty$.
\subsection{H\"older exponent of the measure}
In this subsection, we estimate the H\"older exponent of the measure $\mu$.

Assume $q=ld+h$. Since $\lim\limits_{i\to\infty} \frac{k_{i,j}}{N_{i-1,j}+n_{i,j}}=\alpha_j$ by \eqref{Eq-3}, \eqref{Eq-5.5}, we can choose $q_0$ large enough such that for all $q\geq q_0$,
\begin{align*}
k_{q,h}&\le (\alpha_h+\eta)(N_{q-1,h}+n_{q,h})\\
        &\le (\alpha_h+\eta)(1+\eta)n_{q,h} \quad \text{by \eqref{Eq-eta}}\\
        &\le \frac{\alpha_h+\eta}{1-\eta}n_{q,h}.
\end{align*}
Hence,
\begin{equation}\label{Eq-eta1}
\frac{n_{q,h}}{n_{q,h}+k_{q,h}}
    \ge \frac{1-\eta}{1+\alpha_h}.
\end{equation}

Let
\begin{align*}
s_0 :&=(1-\eta)\Bigl(d-1 + (1-\eta)\min_{1\le j\le d}\bigl\{\dim_{\mathcal H}(E_{\beta_j})\bigr\}\Bigr)\\
&=(1-\eta)\Bigl(d-1 + (1-\eta)\min_{1\le j\le d}\Bigl\{\frac{1}{1+\alpha_j}\Bigr\}\Bigr) < s,
\end{align*}
we will show that for all $q\ge q_0$,
\begin{equation}\label{eq3.31}
\mu(J_q) \lesssim |J_q|^{s_0}.
\end{equation}

Assume that $q=ld+h$ with $0\le h<d$, and that \eqref{eq3.31} holds for $q-1\geq q_0$, that is, $\mu(J_{q-1})\leq |J_{q-1}|^{s_0}\lesssim (\beta_h^{-N_{q-1,h}})^{s_0}$.
We have
\begin{align*}
\frac{-\log_{\beta_h} \mu(J_q)}{-\log_{\beta_h}|J_{q}|}&\geq\frac{-\log_{\beta_h}\mu(J_{q-1})-\log_{\beta_h}(\prod_{j=1}^{d} \frac{1}{\# R_{M,j}^{p_{q,j}}})}{N_{q-1,h}+n_{q,h}+k_{q,h}}\\
&\geq \frac{-\log_{\beta_h}\mu(J_{q-1})+\sum_{j=1}^{d}p_{q,j}M(1-\eta)}{N_{q-1,h}+n_{q,h}+k_{q,h}} \quad\text{by\ \eqref{Eq-M}}\\
&\geq \min\Big\{ \frac{-\log_{\beta_h}\mu(J_{q-1})}{N_{q-1,h}}, \frac{n_{q,h}(1-\eta)+\sum_{j\neq h}n_{q,j}(1-\eta)}{n_{q,h}+k_{q,h}}\Big\}\\
&\geq \min\Big\{ s_0, \frac{n_{q,h}(1-\eta)+(d-1)(n_{q,h}+k_{q,h})(1-\eta)}{n_{q,h}+k_{q,h}}\Big\}\\
&\geq \min\Big\{ s_0, (1-\eta)\Big(\frac{1-\eta}{1+\alpha_h}+d-1\Big)\Big\}\geq s_0.
\end{align*}

Let $\bos n=(n_1,\cdots,n_d)$ with ${\beta_i}^{n_i}\asymp {\beta_j}^{n_j}$ for all $i\neq j$. We aim to estimate the measure of general approximative squares
$$J_{\bos n}:=\prod_{1 \le j \le d} I_{n_{j}}(\varepsilon^{(j)}_{1},\cdots,\varepsilon^{(j)}_{n_j}).$$
Without loss of generality, assume that there exists $q=ld+h(0\leq h<d)$ such that $N_{q-1,h}<n_h<N_{q,h}$.

%We first estimate the length of $I_{n_j}(\varepsilon^{(j)}_{1},\cdots,\varepsilon^{(j)}_{n_j})$
For any $j\neq h$, write $n_j=p_jM+v_j$. Then
\begin{itemize}
\item If $v_j=0$,
$$|I_{n_{j}}(\varepsilon^{(j)}_{1},\cdots,\varepsilon^{(j)}_{n_j})|=|J_{q-1}|\cdot \beta_{j}^{-p_jM},$$
\item If $v_j\neq 0$,
$$|I_{n_{j}}(\varepsilon^{(j)}_{1},\cdots,\varepsilon^{(j)}_{n_j})|\geq \beta_j^{-M} |J_{q-1}|\cdot \beta_{j}^{-(p_j)M},$$
\end{itemize}
Hence, in both cases, we have $|J_{\bos n}|\geq|J_{q-1}|\beta_j^{-(p_j+1)M}$.

\noindent{Case 1:} $N_{q-1,h}<n_h<N_{q-1,h}+n_{q,h}$. \\
Write $n_h=p_hM+v_h$ with $0\leq p_h\leq p_{q,h}$ and $v_h\leq M-1$.

If $v_h= 0$, then
$$|J_{\bos n}|\asymp |J_{q-1}|\beta_h^{-p_hM}.$$
By the definition of the measure $\mu$, we have
\begin{align*}
\mu(J_{\bos n})&=\mu(J_{q-1})\prod^{d}_{j=1}\frac{1}{\# R_{M,j}^{p_{j}}}\lesssim |J_{q-1}|^{s_0}\cdot \prod^{d}_{j=1}\beta_j^{-p_jM(1-\eta)}.
\end{align*}
Thus,
$$\mu(J_{\bos n})\leq |J_{q-1}|^{s_0}\beta_h^{-p_hM(1-\eta)}\cdot\prod_{j\neq h}\beta_{j}^{M}\beta_j^{-(p_j+1)M(1-\eta)} \\
\leq \prod_{j\neq h}\beta_{j}^{(d-1)M}|J_{\bos n}|^{s_0}.$$

If $v_h\neq 0$, we consider the following two cases:\\
\noindent{(i)}When $p_h<p_{q,h}$, we have
$$|J_{\bos n}|\gtrsim \beta_h^{-M}|J_{q-1}|\beta_h^{-p_hM}.$$
By the definition of the measure $\mu$, we have
\begin{align*}
\mu(J_{\bos n})&=\mu(J_{q-1})\prod^{d}_{j=1}\frac{1}{\# R_{M,j}^{p_{j}}}\lesssim |J_{q-1}|^{s_0}\cdot \prod^{d}_{j=1}\beta_j^{-p_jM(1-\eta)}.
\end{align*}
Thus,
$$\mu(J_{\bos n})\lesssim |J_{q-1}|^{s_0}\beta_h^{-p_hM(1-\eta)}\prod_{j\neq h}\beta_{j}^{M}\beta_j^{-(p_h+1)M}
\leq \prod^{d}_{j=1}\beta_{j}^{M}|J_{\bos n}|^{s_0}.$$
\noindent{(ii)}When $p_h=p_{q,h}$, we have
$$|J_{\bos n}|\gtrsim \beta_h^{-M-l_h-1}|J_{q-1}|\beta_h^{-p_hM}.$$
By the definition of the measure $\mu$, we have
\begin{align*}
\mu(J_{\bos n})=\mu(J_{q-1})\prod^{d}_{j=1}\frac{1}{\# R_{M,j}^{p_{j}}}\lesssim |J_{q-1}|^{s_0}\cdot \prod^{d}_{j=1}\beta_j^{-p_jM(1-\eta)}.
\end{align*}
Using $s_0<1-\eta+d-1$, we have
$$\mu(J_{\bos n})\lesssim |J_{q-1}|^{s_0}\beta_h^{-p_hM(1-\eta)}\prod_{j\neq h}\beta_{j}^{M}\beta_j^{-(p_h+1)M}
\leq \beta_h^{l_h+1}\prod^{d}_{j=1}\beta_{j}^{M}|J_{\bos n}|^{s_0}.$$
\noindent{Case 2:} $N_{q-1,h}+n_{q,h}<n_h<N_{q,h}$.\\
We have
$$|J_{\bos n}|\gtrsim |J_{q-1}|\beta_h^{-n_{q,h}-k_{q,h}}.$$
In addition, for $j\neq h$, we also have
$$|J_{\bos n}|\gtrsim |J_{q-1}|\beta_j^{(-p_{j}+1)M}.$$
By the definition of the measure $\mu$, we have
\begin{align*}
\mu(J_{\bos n})=|J_{q-1}|^{s_0}\cdot \beta_h^{-p_{q,h}M(1-\eta)}\prod_{j\neq h}\beta_j^{-p_jM(1-\eta)}\lesssim \beta_h^M|J_{q-1}|^{s_0}(\beta_h^{-n_{q,h}})\prod_{j\neq h}\beta_{j}^{M}\beta_j^{-(p_j+1)M}.
\end{align*}
Thus, by \eqref{Eq-eta1}, we have
\begin{align*}
\mu(J_{\bos n})
&\lesssim \beta_h^M|J_{q-1}|^{s_0}(\beta_h^{-n_{q,h}-k_{q,h}})^{(1-\eta)(\frac{1-\eta}{1+\alpha_h})}\prod_{j\neq h}\beta_{j}^{M}\beta_j^{-(p_j+1)M} \\
&\leq \prod_{j= h}\beta_{j}^{M}|J_{\bos n}|^{s_0}.
\end{align*}
Repeating the above estimations for all coordinate directions, we conclude that
$$\mu(J_{\boldsymbol{n}})\lesssim|J_{\boldsymbol{n}}|^{s_0}.$$

Finally, for any $\bos x=(x_1,\cdots,x_d)\in \mathcal{B}_\infty$, consider $\beta_h^{-n-1}\leq r< \beta_h^{-n}$. The ball $B(\bos x,r)$ intersect at most $\max\limits_{1\leq j\leq d}\{\beta_{j}\}^d$ approximative squares of order $n$. Thus
$$\mu(B(\bos x,r))\lesssim \max\limits_{1\leq j\leq d}\{\beta_{j}\}^d\cdot \beta_h^{-ns_0}\lesssim r^{s_0}.$$
Applying Proposition \ref{Mdp}, we conclude that $\dim_\mathcal{H}(E_{\beta_j})\geq s_0$.
By letting $\eta\to 0$, we complete the proof of the desired lower bound.
\section*{Acknowledgment}
This work is supported by National Natural Science Foundation of China (No.12371086, 12271175) and National Key R\&D Program of China (No. 2024YFA1013701).

\end{document}